\documentclass[final,10pt,a4paper,leqno]{amsart}
\usepackage[english]{babel}
\usepackage[T1]{fontenc}
\usepackage{amssymb,cite}
\usepackage{amsfonts}
\usepackage{amsmath}
\usepackage{mathptm}
\usepackage[notref,notcite]{showkeys}
\usepackage[dvips]{graphicx}
\renewcommand{\rmdefault}{ptm}
\renewcommand{\ttdefault}{pcr}
\renewcommand{\sfdefault}{phv}
\newcommand{\R}{\mathbb{R}}
\newcommand{\zbar}{\bar z}
\newcommand{\ol}{\overline}
\newcommand{\loc}{\mathrm{loc}}
\usepackage{eucal}
\begin{document}

\thispagestyle{empty}
\date{}
\title{{{Mean value type inequalities for quasinearly\\ subharmonic functions}}}
\author{Oleksiy Dovgoshey}
\address{Institute of Applied Mathematics and Mechanics of NASU\\
R. Luxemburg str.~74\\
Donetsk 83114, Ukraine} \email{aleksdov@mail.ru}
\author{Juhani Riihentaus}
\address{Department of Physics and Mathematics\\
University of Joensuu\\
P.O. Box 111\\
FI-80101 Joensuu, Finland} \email{juhani.riihentaus@gmail.com}
\begin{abstract}
The mean value inequality is characteristic for upper
semicontinuous functions to be subharmonic. Quasinearly
subharmonic functions generalize subharmonic functions.  We find
the necessary and sufficient conditions under which subsets of
balls are big enough for the catheterization of nonnegative,
quasinearly subharmonic functions by mean value inequalities.
Similar result is obtained also for generalized mean value
inequalities where, instead of balls, we consider arbitrary
bounded sets which have nonvoid interiors and instead of the
volume of ball some functions depending on the radius of this
ball.
\end{abstract}
 \maketitle

\section{Subharmonic functions and generalizations. Some definitions and results}
\subsection{Notation.}
Our notation is rather standard, see e.g.
\cite{Ri07$_1$,Ri07$_2$,PaRi08,Ri08} and the references therein.
 If $E\subset {\mathbb{R}}^n$
and $x\in {\mathbb{R}}^n$, then we write $\delta_E(x):=\inf
\{\,\mid x-y\mid \,:\, y\in E^c \,\}$, where
$E^c={\mathbb{R}}^n\setminus E$. The Lebesgue measure in
${\mathbb{R}}^n$ is denoted by $m_n$. We write $B^n(x,r)$ for the
open ball in ${\mathbb{R}}^n$ with center $x$ and radius $r$.
Recall that  $m_n(B^n(x,r))=\nu_nr^n$, where
$\nu_n:=m_n(B^n(0,1))$. We denote by $Int\, D,\ \overline D$ and
$\partial D$ the interior, the closure and the boundary of a set
$D\subseteq\mathbb R^n$, i.e., $\overline{B^n(x,r)}$ is the closed
ball with center $x$ and radius $r$. Note also that our constants
$C$ and $K$ are nonnegative, mostly $\geq 1$, and may vary from
line to line.
\subsection{Subharmonic functions and generalizations}
Let $\Omega$ be an open set in ${\mathbb{R}}^n$, $n\geq 2$. Let $u:\Omega \rightarrow [-\infty ,+\infty )$ be Lebesgue measurable. We adopt the following definitions:
\begin{itemize}
\item[(i)] $u$ is \emph{subharmonic} if $u$ is upper
semicontinuous and if
\begin{equation} \label{e1*}
u(x)\leq \dfrac{1}{\nu_nr^n}\int\limits_{B^n(x,r)}u(y)\,
dm_n(y)
\end{equation}
for all balls $\overline{B^n(x,r)}\subset \Omega$.  A subharmonic
function may be $\equiv -\infty$ on any component of $\Omega$. See
\cite[p.~1]{Ra37},  \cite[p.~18]{Br69}, \cite[p.~9]{He71},   and
\cite[p.~60]{ArGa01}. \item[(ii)] $u$ is \emph{nearly subharmonic}
if $u^+\in{\mathcal{L}}^1_{{\textrm{loc}}}(\Omega)$ and inequality
\eqref{e1*} holds for all balls $\overline{B^n(x,r)}\subset
\Omega$. Observe that this definition is slightly more general
than the usual one, compare \cite[p.~51]{Ri07$_2$},  with the
standard definitions \cite[p.~20]{Ra37}, \cite[p.~30]{Br69},  and
\cite[p.~14]{He71}. \item[(iii)] Let $K\geq 1$. Then $u$ is
\emph{$K$-quasinearly subharmonic} if $u^+\in
{\mathcal{L}}^1_{{\textrm{loc}}}(\Omega)$ and inequality
\begin{equation*}
u_M(x)\leq \dfrac{K}{\nu_nr^n}\int\limits_{B^n(x,r)}u_M(y)\,
dm_n(y)\end{equation*} holds for all $M\geq 0$ and for all balls
$\overline{B^n(x,r)}\subset \Omega$. Here $u_M:=\max \{\,
u,-M\,\}+M$.

The function $u$ is \emph{quasinearly subharmonic} if $u$ is
$K$-quasinearly subharmonic for some $K\geq 1$. For the definition
and properties of quasinearly subharmonic functions, see e.g.
\cite{Ri89,Pa94,Ri00,Ri01,Ri06,Ri07$_1$,Ri07$_2$,Ko07,
DjPa07,PaRi08,Ri08},  and the references therein. We write
$QNS(\Omega)$ for the set of all nonnegative quasinearly
subharmonic functions on the open set $\Omega\subseteq\mathbb
R^n$.
\end{itemize}
\subsection{}\label{p1.3} {\textbf{Proposition.}} (cf. \cite{Ri07$_2$}, Proposition~2.1, pp.~54-55)
\emph{The following statements hold:}
\begin{itemize}
\item[(i)] \emph{A subharmonic function is nearly subharmonic but
not conversely.} \item[(ii)] \emph{A function is nearly
subharmonic if and only if it is $1$-quasinearly subharmonic.}
\item[(iii)] \emph{A nearly subharmonic function is quasinearly
subharmonic but not conversely.}
\item[(iv)] \emph{If $u:
\Omega\rightarrow [0,+\infty )$ is Lebesgue measurable, then $u$
is $K$-quasinearly subharmonic if and only if}
$u\in{\mathcal{L}}^1_{{\textrm{loc}}}(\Omega)$ \emph{and
\begin{equation}\label{e1}
u(x)\leq \frac K{m_n(B^n(x,r))}\int\limits_{B^n(x,r)}u(y)\,dm_n(y)
\end{equation}
for all closed balls $\overline{B^n(x,r)}\subset \Omega$.}
\end{itemize}
Note that if $u$ is $K$-quasinearly subharmonic and nonnegative in
$\Omega$, then \eqref{e1} holds also for every open ball
$B^n(x,r)\subseteq\Omega$.

Let $A$ be a subset of the open half-line $(0,\infty)$ such that
$0$ is a limit point of $A$ and let $u:\Omega\to[-\infty,+\infty)$
be an upper-semicontinuous function on an open set
$\Omega\subseteq\mathbb R^n$. The classical Blascke--Privalov
theorem, see, for example, \cite[Chapter II,\S2]{Br} implies that
$u$ is subharmonic if  inequality \eqref{e1*} holds whenever $r\in
A$ and $\overline{B^n(x,r)}\subset\Omega$. Moreover, the simple
examples show that if nonnegative $u\in\mathcal
L^1_{\loc}(\Omega)$, then the fulfilment of \eqref{e1} for all
$(x,r)\in\Omega\times A$ with $\overline{B^n(x,r)}\subset\Omega$
does not, generally, imply $u\in QNS(\Omega)$. A legitimate
question to raise in this point is in finding the sets
$A\subseteq(0,\infty)$ for which every nonnegative  $u\in \mathcal
L^1_{\loc}(\Omega)$ is quasinearly subharmonic if \eqref{e1} holds
for   $(x,r)\in\Omega\times A$ whenever
$\overline{B^n(x,r)}\subset\Omega$.

\subsection{Definition}\label{d3.5}{\it Let $\Omega$ be an open set in $\mathbb R^n$.
A set $A\subseteq(0,\infty)$ is favorable for $\Omega$ (favorable
for the characterization of nonnegative, quasinearly subharmonic
functions in $\Omega$) if for every nonnegative $u\in\mathcal
L_{\rm loc}^1(\Omega)$ the following conditions are equivalent:
\begin{itemize}
\item[(i)] $u\in QNS(\Omega)$;

\item[(ii)] There is $K=K(u,A,\Omega)\geq1$ such that for all
$x\in\Omega$ the inequality
\begin{equation}\label{e8}
u(x)\leq\frac K{\nu_nr^n}\int\limits_{B^n(x,r)}u(y)\,dm_n(y)
\end{equation}
holds whenever $r\in A$ and $\overline{B^n(x,r)}\subset\Omega$.
\end{itemize}
}

We can characterize the favorable subsets of $(0,\infty)$ by the
following way.

\subsection{Theorem}\label{t1.5}
{\it The following three statements are equivalent for every
$A\subseteq(0,\infty)$:
\begin{itemize}

\item[(i)] $A$ is favorable for all open sets
$\Omega\subseteq\mathbb R^n$;

\item[(ii)]  The characteristic function
$$
\chi_\Gamma(x)=\begin{cases} 1&\text{if }x\in\Gamma\\
0&\text{if } x\in\Omega\setminus\Gamma
\end{cases}
$$
is quasinearly subharmonic for all open sets
$\Omega\subseteq\mathbb R^n$ and all Lebesgue measurable sets
$\Gamma\subseteq\Omega$ if  and only if there is a constant
$K=K(\Gamma,\Omega,n)$ such that the inequality
\begin{equation}\label{e3}
m_n(B^n(x,r))\leq Km_n(\Gamma\cap B^n(x,r))
\end{equation}
holds for $(x,r)\in\Omega\times A$ whenever
$\overline{B^n(x,r)}\subset\Omega$.

\item[(iii)] There exists $C=C(A)>1$ such that
$$
\Big[\frac xC,x\Big]\cap A\ne\emptyset
$$
for every $x\in(0,\infty)$.
\end{itemize}}

We shall prove the equivalence (i)$\equiv$(iii) in Theorem
\ref{t3.6} below. Observe also that the implication
(i)$\Rightarrow$(ii) is trivial and that (ii)$\Rightarrow$(iii)
follows directly from the proof of Theorem \ref{t3.6}.  The
quasidisks give the important example of the sets $\Gamma$ such
that \eqref{e3} holds in a bounded domain $\Omega\subseteq\mathbb
R^2$ whenever $\overline{B^2(x,r)}\subset\Omega$. It is a
particular case of the Gehring--Martio result which proves
\eqref{e3} for the so-called quasiextremal distance domains in
$\mathbb R^n,\ n\geq2$. See \cite[Lemma~2.13]{GM}.

The following  result closely connected to Theorem \ref{t1.5},
follows from Theorem \ref{p2.13} formulated in the second section
of the paper.

\subsection{Theorem}\label{t1.6}
{\it Let $f$ be a positive function on $(0,\infty)$. The following
three statements are equivalent. \begin{itemize}

\item[\rm(i)] For all open sets $\Omega\subseteq\mathbb R^n$,
Lebesgue measurable functions $u:\Omega\to[0,\infty)$ are
quasinearly subharmonic if and only if there are constants
$K=K(u,\Omega,n)\geq1$ such that
$$
u(x)\leq \frac K{(f(r))^n}\int\limits_{B^n(x,r)}u(y)\,dm_n(y)
$$
for all closed balls $\overline{B^n(x,r)}\subset\Omega$.

\item[\rm(ii)] For all open sets $\Omega\subseteq\mathbb R^n$ and
all Lebesgue measurable sets $\Gamma\subseteq\Omega$ the
characteristic functions $ \chi_\Gamma$ are quasinearly
subharmonic if and only if there are constants
$K=K(\Gamma,\Omega,n)$ such that the inequality
$$
(f(r))^n\leq K m_n({B^n(x,r)\cap\Gamma})
$$
holds for all closed balls $\overline{B_n(x,r)}\subset\Omega$ with
$x\in\Gamma$.

\item[\rm(iii)] There are a set $A\subseteq(0,\infty)$ and a
constant $c>1$ such that:
\begin{itemize}

\item[\rm(iii$_1$)] The inequality $f(r)\leq cr$ holds for all
$r\in(0,\infty)$;

\item[\rm(iii$_2$)] $\ln A$ is an $\varepsilon$-net in $\mathbb R$
for some $\varepsilon>0$;

\item[\rm(iii$_3$)] The inequality
$$
\frac1cr\leq f(r)
$$
holds for all $r\in A$.
\end{itemize}
\end{itemize}}

Note that condition (iii) of Theorem \ref{t1.5} holds if and only
if the set $\ln(A):=\{\ln x:x\in A\}$ is an $\varepsilon$-net in
$\mathbb R$ for some $\varepsilon>0$. A  characterization in terms
of porosity for the sets $A$ which are favorable for {\it bounded}
open sets $\Omega\subseteq\mathbb R^n$ is proved in
Theorem~\ref{t2.11} below.

\section{Generalized mean value inequalities}
Inequality \eqref{e1}, characteristic  for quasinearly subharmonic
functions,  can be  generalized by some distinct ways. Our first
theorem characterizes nonnegative quasinearly subharmonic
functions via mean values over some sets more general than just
balls.

\subsection{Similarities of the Euclidean space}
Let $\Omega$  and $D$ be  subsets of $\mathbb R^n$ with marked
points $p_\Omega\in\Omega$ and $p_D\in D$. In what follows we
always suppose that $Int\, D\ne\emptyset$ and $p_D\in Int\,D$.
Denote by $Sim(p_D,p_\Omega)$ the set of all similarities
$h:\mathbb R^n\to\mathbb R^n$ such that $h(p_D)=p_\Omega$ and
$h(D)\subseteq\Omega$. Recall that $h$ is a similarity if there is
a positive number $k=k(h)$, the similarity constant of $h$, such
that
$$
|h(x)-h(y)|=k|x-y|
$$
for all $x,y\in\mathbb R^n$. The group of all similarities of the
Euclidean space $\mathbb R^n$ is sometimes denoted as $SM(\mathbb
R^n)$, see e.g. \cite[5.1.14]{JD}, and we also adopt this
designation. Observe that each similarity $h\in SM(\mathbb R^n)$
can be written in the form
$$
h(x)=k(h)Tx+a,\quad x\in\mathbb R^n,
$$
where $k(h)>0,$ and $T:\mathbb R^n\to\mathbb R^n$ is an orthogonal
linear mapping and $a\in\mathbb R^n$.

\subsection{Theorem}\label{t2.2}
\emph{Let $\Omega$ be an open set in $\mathbb R^n,\ n\geq2$, let
$D$ be a bounded, Lebesgue measurable set with the marked point
$p_D\in Int\,D$ and let $u:\Omega\to[0,\infty)$ be a function from
$\mathcal L_{\rm loc}^1(\Omega)$. Then $u$ is quasinearly
subharmonic if and only if there is $C\geq1$ such that
\begin{equation}\label{e2} u(x_\Omega)\leq\frac
C{m_n(h(D))}\int\limits_{h(D)}u(y)\,dm_n(y)
\end{equation} for every point $x_\Omega$ and all $h\in Sim(p_D,x_\Omega)$.
If $u$ is $K$-quasinearly subharmonic, then $C=C(D,p_D,K,n)$ and,
conversely, if \eqref{e2} holds, then $u$ is $K$-quasinearly
subharmonic with $K=K(D,p_D,C,n)$.}

\begin{proof}
Write
\begin{equation}\label{2*}
R_D:=\sup_{y\in D}|p_D-y|,\quad\text{and}\quad
r_D:=\delta_{Int(D)}(p_D).
\end{equation}

Suppose that $u$ is quasinearly subharmonic, i.e., there is $K
\geq1$ such that \eqref{e1} holds for all
$B^n(x,r)\subseteq\Omega$. Let $x_\Omega$ be an arbitrary point of
$\Omega$ and let $h\in Sim(p_D,x_\Omega)$. The last membership
relation implies the inclusions
$$
B^n(x_\Omega,k(h)r_D)\subseteq \Omega \quad\text{and}\quad
h(D)\subseteq\overline{ B^n(x_\Omega,k(h)R_D)}
$$
where $k(h)$ is the similarity constant of $h$. Consequently we
obtain
\begin{multline*}
\frac1{m_n(h(D))}\int\limits_{h(D)}u(y)\,dm_n(y)\geq
\frac1{\nu_n(k(h)R_D)^n}\int\limits_{h(D)}u(y)\,dm_n(y)\\
\geq\Big(\frac{r_D}{R_D}\Big)^n
\frac1{\nu_n(k(h)r_D)^n}\int\limits_{B^n(x_\Omega,k(h)r_D)}u(y)\,dm_n(y)\geq
\Big(\frac{r_D}{R_D}\Big)^n\frac{u(x_\Omega)}K.
\end{multline*}
Thus if $f$ is $K$-quasinearly subharmonic, then \eqref{e2} holds
with
$$C=\frac{K(R_D)^n}{(r_D)^n}.$$

Conversely, suppose that \eqref{e2} holds with some $C\geq1$ for
all $x_\Omega$ and all $h\in Sim(p_D,x_\Omega)$. Let
$B^n(x_\Omega,r_0)\subseteq\Omega$. Let $h$ be an arbitrary
similarity with  $k(h)=\frac{r_0}{R_D}$ and with
$h(p_D)=x_\Omega$. Then we have $h\in Sim(p_D,x_\Omega)$ and
$B^n(x_\Omega,k(h)r_D)\subseteq h(D)\subseteq\overline{
B^n(x_\Omega,r_0)}$. Consequently
$$
\frac
C{m_n(B^n(x_\Omega,r_0))}\int\limits_{B^n(x_\Omega,r_0)}u(y)\,dm_n(y)\geq
\frac{Cm_n(h(D))}{m_n(B^n(x_\Omega,r_0))m_n(h(D))}\int\limits_{h(D)}u(y)\,dm_n(y)\geq
u(x_\Omega)\frac{m_n(h(D))}{m_n(B^n(x_\Omega,r_0))}.
$$
Since
$$
\frac{m_n(B^n(x_\Omega,r_0))}{m_n(h(D))}=\frac{\nu_n(r_0)^n}{m_n(h(D))}=
\frac{\nu_n(k(h))^n(R_D)^n}{(k(h))^nm_n(D)}=\frac{\nu_nR_D^n}{m_n(D)},
$$
inequality (1) holds with
$$
K=C\frac{\nu_nR_D^n}{m_n(D)}.
$$
\end{proof}
\subsubsection{Remark}\label{r2.2.1}
The standard notion of quasinearly subharmonicity is defined by
the condition \eqref{e2}, where $D=B^n(0,1)$, $p_D=0$ and the
considered similarities $h$ are of the form $h(x)=r_0x+x_\Omega$.
The point of Theorem \ref{t2.2} is that the definition and its
consequences are, however, much more general: Instead of just
$D=B^n(0,1)$ and $p_D=0$ one may consider arbitrary bounded sets
$D$ with nonvoid interior $Int\, D$.

\subsubsection{Remark}\label{r2.2.2}
Inequality \eqref{e2} remains valid for each nonnegative
quasinearly subharmonic function if we use bi-Lipschitz mappings
$h$ instead of similarities, but in this more general case the
constant in \eqref{e2} depends on the Lipschitz constant of $h$.
See Lemma~2.1 in \cite{DR}.
\medskip

Inequality \eqref{e2} remains also valid  for unbounded sets $D$
if $m_n(D)<\infty$.
\subsection{Proposition}\label{p2.3}
{\it Let $\Omega$ be an open set in $\mathbb R^n,\ n\geq2$, $D$ a
Lebesgue measurable set with $m_n(D)<\infty$, $p_D$ a  point of
$Int(D)$ and let $u:\Omega\to[0,\infty)$ be a $K$-quasinearly
subharmonic function. Then there is a constant $C=C(D,p_D,K,n)$
such that \eqref{e2} holds for all $x_\Omega\in\Omega$ and $h\in
Sim(p_D,x_\Omega)$.}
\begin{proof}
If $D$ is bounded, then this proposition follows from Theorem~2.2.
Suppose $D$ is an unbounded. Let $t>1$ be a constant. It is easy
to show that there is a ball $B^n(p_D,r_t)$ with a sufficiently
large radius $r_t$ such that
\begin{equation}\label{e4}
tm_n(D\cap B^n(p_D,r_t))\geq m_n(D).
\end{equation}
Write
$$
D_t:=D\cap B^n(p_D,r_t)\quad\text{and}\quad p_{D_t}:=p_D.
$$
Note that $D_t$ satisfies all conditions of Theorem 2.2 and that
$p_{D_t}\in Int(D_t)$. Consequently there is $K\geq1$ such that
the inequality
$$
u(x_\Omega)\leq\frac
K{m_n(h(D_t))}\int\limits_{h(D_t)}u(y)\,dm_n(y)
$$
holds for all $x_\Omega$ and $h\in Sim(p_{D_t},x_\Omega)$. If
$h\in Sim(p_{D},x_\Omega)$, then we have $h\in
Sim(p_{D_t},x_\Omega)$ and $h(D_t)\subseteq h(D)$. Since
$$
\frac{m_n(D_t)}{m_n(D)}=\frac{m_n(h(D_t))}{m_n(h(D))}
$$
for all $h\in SM(\mathbb R^n)$, \eqref{e4} implies the inequality
$$
\frac1{m_n(h(D_t))}\int\limits_{h(D_t)}u(y)\,dm_n(y)\leq \frac
t{m_n(h(D))}\int\limits_{h(D)}u(y)\,dm_n(y).
$$
Thus \eqref{e2} holds for all $h\in Sim(p_{D},x_\Omega)$ with
$C=tK$.
\end{proof}
\subsubsection{Remark}\label{r2.3.1}
If $ Sim(p_{D},x_\Omega)=\emptyset$ for all $x_\Omega$, then
Proposition \ref{p2.3} is vacuously true.
\medskip

Let $ \varphi: SM(\mathbb R^n)\to(0,\infty) $ be a function such
that the equality
$$
\varphi(h)=\varphi( is\circ h)
$$
holds for all $h\in SM(\mathbb R^n)$ and for all isometries
$is:\mathbb R^n\to\mathbb R^n$. Then we have
$\varphi(h_1)=\varphi(h_2)$ whenever $k(h_1)=k(h_2)$, that is
there is  a function $f:(0,\infty)\to(0,\infty)$ such that the
equality
\begin{equation}\label{e5}
\varphi(h)=f(k(h))
\end{equation}
is fulfilled for all $h\in SM(\mathbb R^n)$ with $k(h)$ equals the
similarity constant of $h$. For instance, if $D$ is a bounded
nonvoid subset of $\mathbb R^n$, we can put
$\varphi(h)=diam(h(D))$. Other  examples can be found in
\eqref{e2.15}--\eqref{e2.17.4}.

Let $D$ be a measurable subset of $\mathbb R^n$ with a marked
point $p_D\in Int\,D$.  For every open set $\Omega\subseteq\mathbb
R^n$ define a subset $Q(f,D,\Omega)\subseteq\mathcal L_{\rm
loc}^1(\Omega)$ by the rule:

{\it $ u\in Q(f,D,\Omega)$ if and only if $u\geq0$ and
$u\in\mathcal L_{\rm loc}^1(\Omega) $ and there is $K=K(u)\geq1$
such that the inequality
\begin{equation}\label{e6}
u(x_\Omega)\leq \frac
K{(f(k(h)))^{n}}\int\limits_{h(D)}u(y)\,dm_n(y)
\end{equation}
holds for every $x_\Omega\in\Omega$ and all $h\in
Sim(p_D,x_\Omega)$.}

It is clear that $QNS(\Omega)=Q(f,D,\Omega)$ if $D$ satisfies the
conditions of Theorem \ref{t2.2}  and $f(k(h))=\break
k(h)(m_n(D))^{\frac1n}$.
\subsection{Proposition}\label{p2.4}
{\it Let $D$ be a bounded, Lebesgue measurable subset of $\mathbb
R^n$  with a marked point $p_D\in Int\,D$ and let
$\varphi:SM(\mathbb R^n)\to(0,\infty)$,
$f:(0,\infty)\to(0,\infty)$  be functions such that \eqref{e5}
holds for all $h\in SM(\mathbb R^n)$. Then the inclusion
\begin{equation}\label{e7}
QNS(\Omega)\subseteq Q(f,D,\Omega)
\end{equation}
is valid for all open sets $\Omega\subseteq\mathbb R^n$ if and
only if there is $c\geq1$ such that  the inequality
\begin{equation}\label{e8*}
f(k)\leq ck
\end{equation}
holds for all $k\in(0,\infty)$.}
\begin{proof}
Suppose that inclusion \eqref{e7} holds for all open sets
$\Omega\subseteq\mathbb R^n$. Let $\Omega$ be an open half-space
of $\mathbb R^n$. Then for every $k_0\in(0,\infty)$ there is a
similarity $h_0$ with the similarity constant $k(h_0)=k_0$ such
that $h_0(D)\subseteq\Omega$. The constant function $u_1,\
u_1(x)\equiv1$ for $x\in\Omega$, belongs to $QNS(\Omega)$. Hence,
by \eqref{e7}, $u_1\in Q(f,D,\Omega)$ and it follows from
\eqref{e6} that
$$
1=u_1(h_0(p_D))\leq\frac
K{(f(k_0))^{n}}\int\limits_{h_0(D)}u_1(x)\,dm_n(x)=
\frac{Km_n(h_0(D))}{(f(k_0))^{n}}=\frac{K(k_0)^nm_n(D)}{(f(k_0))^{n}}.
$$
Consequently \eqref{e7} implies \eqref{e8*} for all
$k\in(0,\infty)$ with
$$
c=(K(u_1)m_n(D))^{\frac1n}\vee1.
$$

Conversely suppose that \eqref{e8*} holds for all
$k\in(0,\infty)$. Then using Theorem 2.2 we obtain the following
inequalities for every open set $\Omega\subseteq\mathbb R^n$,
every $u\in QNS(\Omega)$, every $x_\Omega\in\Omega$ and every
$h\in Sim(p_D,x_\Omega)$:
\begin{multline*}
u(x_\Omega)\leq\frac{C(u)}{m_n(h(D))}\int\limits_{h(D)}u(x)\,dm_n(x)=
\frac{C(u)(f(k(h)))^{n}}{(k(h))^nm_n(D)(f(k(h)))^{n}}\int\limits_
{h(D)}u(x)\,dm_n(x)\\\leq
\frac{C(u)c^n}{m_n(D)(f(k(h)))^{n}}\int\limits_{h(D)}u(x)\,dm_n(x).
\end{multline*}
Hence \eqref{e6} holds with
$$
K=\frac{C(u)c^n}{m_n(D)}\vee1.
$$
Thus \eqref{e7} is valid for all open sets $\Omega\subseteq\mathbb
R^n$.
\end{proof}
Before passing to the equality
$$
Q(f,D,\Omega)= QNS(\Omega)
$$
we consider one relevant question.

\subsection{Theorem}\label{t3.6}{\it Let $A$ be a subset of $(0,\infty)$.
Then $A$ is favorable for all open sets $\Omega\subseteq\mathbb
R^n$ if and only if the following statement holds.
\begin{itemize}
\item[(s)] There exists $C=C(A)>1$ such that
\begin{equation}\label{e9}
\Big[\frac xC,x\Big]\cap A\ne\emptyset
\end{equation}
for every $x\in(0,\infty)$. \end{itemize}}\medskip

The following lemma will be used in the proof of Theorem
\ref{t3.6}.

\subsection{Lemma}\label{l3.8}{\it Let $A\subseteq(0,\infty)$.
Statement (s) of Theorem \ref{t3.6} does not hold with this $A$,
if and only if there are disjoint open intervals $(a_m,b_m),\
a_m<b_m,\ m=1,2,\dots$, in $(0,\infty)\setminus A$ such that
\begin{equation}\label{eq12}
\lim_{m\to\infty}\frac{a_m}{b_m}=0
\end{equation}
and either
\begin{equation}\label{eq13}
\lim_{m\to\infty}a_m=\lim_{m\to\infty}b_m=0
\end{equation}
or
\begin{equation}\label{eq14}
\lim_{m\to\infty}a_m=\lim_{m\to\infty}b_m=\infty.
\end{equation}}

\begin{proof}
If statement (s) holds, then using \eqref{e9} we obtain that
$$
\frac ab\geq\frac1{C(A)}
$$
for every open interval $(a,b)$ in $(0,\infty)\setminus A$. This
inequality contradicts \eqref{eq12}.

Conversely, suppose that statement (s) of Theorem \ref{t3.6} does
not hold and that $0$ and $\infty$ are limit points of $A$. Then
for every natural $i\geq2$ there is $x\in(0,\infty)$ such that
$$
\Big(\frac xi,x\Big)\cap A=\emptyset.
$$
Let $\bar A$ be the closure of $A$ in $(0,\infty)$. Write
$(a_i,b_i)$ for the connected component of
$(0,\infty)\setminus\bar A$ which contains $(\frac xi,x)$.
 Since both $0$ and $\infty$
are the limit points of $A$ we have
$$
0<a_i<b_i<\infty.
$$
Passing to convergent, in $[0,\infty]$, subsequences
$\{a_{i_m}\}_{m\in\mathbb N}$ and $\{b_{i_m}\}_{m\in\mathbb N}$ it
is easy to see that limits $\lim_{m\to\infty}a_{i_m}$ and
$\lim_{m\to\infty}b_{i_m}$ are $0$ or $\infty$ and that the
equalities
$$
\lim_{m\to\infty}a_{i_m}=0\quad\text{and}\quad
\lim_{m\to\infty}b_{i_m}=\infty.
$$
cannot be true simultaneously. Renaming $a_m:=a_{i_m}$ and
$b_m:=b_{i_m}$ we obtain the desirable sequence of intervals in
$(0,\infty)\setminus  A$.

If at least one of the points $0$ and $\infty$ is not a limit
point of $A$, then there is $\varepsilon>0$ such that
$$
A\subset(0,\varepsilon]\quad\text{or}\quad
A\subset[\varepsilon,\infty).
$$
Each of these inclusions implies evidently the existence of
desired intervals in $(0,\infty)\setminus\bar A$.
\end{proof}
\begin{proof}[Proof of Theorem \ref{t3.6}]
We shall first prove that $A$ is favorable for all open sets
$\Omega\subseteq\mathbb R^n$ if statement (s) holds.

Suppose that (s) is true. Let $\Omega$ be an open set in $\mathbb
R^n$ and let $u\in\mathcal L_{\rm loc}^1(\Omega)$ be  a
nonnegative function which satisfies condition (ii) of Definition
\ref{d3.5}. It is enough to show that $u\in QNS(\Omega)$. To prove
this, consider an arbitrary
$\overline{B^n(x_\Omega,r_0)}\subseteq\Omega$. By statement (s)
there is $r_1\in A$ such that
$$
\frac{r_0}C\leq r_1\leq r_0
$$
where the constant $C=C(A)>1$. Using this double inequality and
condition (ii) of Definition \ref{d3.5} we obtain
$$
u(x_\Omega)\leq\frac
K{\nu_n(r_1)^n}\int\limits_{B^n(x_\Omega,r_1)}u(y)\,dm_n(y)\leq
\frac{KC^n}{\nu_n(r_0)^n}\int\limits_{B^n(x_\Omega,r_0)}u(y)\,dm_n(y).
$$
Statement (iv) of Proposition 1.3 implies that $u\in QNS(\Omega)$.

Conversely, suppose that $A$ is favorable for every open set
$\Omega\subseteq\mathbb R^n$. We must show that  (s)  holds. If
(s) does not hold then, by Lemma~\ref{l3.8}  there is a sequence
of disjoint open intervals in $(0,\infty)$ satisfying \eqref{eq12}
and \eqref{eq13} or \eqref{eq12} and \eqref{eq14}. Suppose that
\eqref{eq12} and \eqref{eq13} hold. Then for every integer $N_0>2$
there is a sequence of  open intervals $(a_m,b_m)$ such that
\begin{equation}\label{e10}
0<b_{m+1}<a_m<2a_m<\frac1{N_0}b_m<b_m
\end{equation}
and
\begin{equation}\label{e11}
(a_m,b_m)\cap A=\emptyset
\end{equation}
for $m=1,2,\dots$ and
\begin{equation}\label{e12}
\lim_{m\to\infty}\frac{b_m}{a_m}=\infty.
\end{equation}
Moreover, passing, if necessary, to a subsequence we may assume
that
\begin{equation}\label{e13}
\sum_{m=1}^\infty b_m<\infty.
\end{equation}
For the sake of simplicity, we shall describe our constructions
only on the plane   but in such a way that a generalization to the
dimensions $n\geq3$ is a trivial matter.

Define the points $z_m\in\mathbb C,\ m=1,2,\dots$, as
$$
z_m:=\begin{cases} 0&\text{if }m=1\\
2\sum_{i=1}^{m-1}b_i&\text{if }m\geq2
\end{cases}
$$
and write
$$
R_1:=\{z\in\mathbb C:0< Re(z)<2b_1,\ |Im(z)|<a_2\}
$$
and
$$
R_m:=\{z\in\mathbb C:2\sum_{i=1}^{m-1}b_i<
Re(z)<2\sum_{i=1}^mb_i,\ |Im(z)|<a_{m+1}\}
$$
for $m\geq2$. Using \eqref{e10} we see that $B^2(z_m,b_m)$ are
open, pairwise disjoint balls and that $R_m$ are open, pairwise
disjoint rectangles. The desired domain $\Omega$ is, by
definition, the union
$$
\bigcup_{m=1}^\infty\Big(B^2\Big(z_m,\frac{b_m}{N_0}\Big)\cup
R_m\Big),
$$
see Fig. 1. Let us define now a function $u$ as the characteristic
function of the set
\begin{equation}\label{13*}
X:=\bigcup_{m=1}^\infty\overline{B^2(z_m,a_m)},
\end{equation}
 i.e.,
\begin{equation}\label{e14}
u(z):=\begin{cases} 1&\text{if }z\in X\\
0&\text{if }z\in\Omega\setminus X.
\end{cases}
\end{equation}

\begin{figure}[htb]
\includegraphics[width=\textwidth,keepaspectratio]{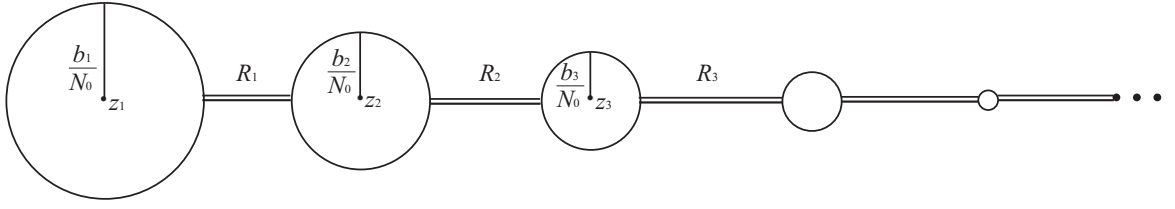}
\caption{Domain $\Omega$ is an infinite sequence of balls
$B^2(z_m,\frac{b_m}{N_0})$ united by the thin rectangles $R_m$.}
\end{figure}
\noindent It is clear that $u\geq0$ and that $u\in\mathcal
L^1(\Omega)$. Moreover, since
\begin{equation}\label{14*}
\frac1{m_2(B^2(z_m,\frac{b_m}{2N_0}))}\int\limits_{B^2(z_m,\frac{b_m}{2N_0})}
u(z)\,dm_2(z)=\frac{4N_0^2(a_m)^2}{(b_m)^2},
\end{equation}
statement (iv) of Proposition 1.3 and limit relation \eqref{e12}
imply $u\not\in QNS(\Omega)$. It remains to show that there is $K$
such that \eqref{e8} holds  whenever $r\in A$ and
$\overline{B^2(x,r)}\subseteq\Omega$. If $x\in\Omega\setminus X$,
then \eqref{e8} is trivial and we must consider only  $x\in X$.
The last membership relation implies that there exists  $m=m_x$
such that
\begin{equation}\label{e15}
x\in\overline{B^2(z_m,a_m)}.
\end{equation}
Let us consider all  $r\in A$ such that
\begin{equation}\label{e16}
\overline{B^2(x,r)}\subseteq\Omega.
\end{equation}
From \eqref{e11} follows that either $r\geq b_{m_x}$ or $r\leq
a_{m_x}$. If $r\geq b_{m_x}$, then we have
\begin{equation}\label{e17}
B^2(x,r)\supseteq\overline{B^2(z_{m_x},\frac{b_{m_x}}{N_0})}.
\end{equation}
Indeed, the triangle inequality and \eqref{e10} imply
$$
|y-x|\leq |x-z_{m_x}|+|z_{m_x}-y|\leq
a_{m_x}+\frac{b_{m_x}}{N_0}<\frac34b_{m_x}<r
$$
for all $y\in \overline{B^2(z_{m_x},\frac{b_{m_x}}{N_0})}$.
Inclusion \eqref{e17} and definition of $\Omega$ show that
$B^2(x,r)\not\subseteq\Omega$ if $r\geq b_{m_x}$. Consequently if
\eqref{e16} holds, then
\begin{equation}\label{e18}
r\leq a_{m_x}.
\end{equation}
Using the last inequality, \eqref{e15} and \eqref{e10} we obtain
$$
\overline{B^2(x,r)}\subseteq
B^2\Big(z_m,\frac{b_m}{N_0}\Big),\quad m=m_x,
$$
for these $x$ and $r$. Hence the equality
\begin{equation}\label{e18*}
\frac1{m_2(B^2(x,r))}\int\limits_{B^2(x,r)}u(y)\,dm_2(y)=\frac{m_2
(B^2(z_m,a_m)\cap B^2(x,r))}{m_2(B^2(x,r))}
\end{equation}
holds for such $x$ and $r$. Write
\begin{equation}\label{e19}
C=\inf\frac{m_2(B^2(z_m,a_m)\cap B^2(x,r))}{m_2(B^2(x,r))}
\end{equation}
where the infimum is  taken over the set of all balls $B^2(x,r)$
with $x\in B^2(z_m,a_m)$ and with $r\leq a_m$. If $r$ is fixed and
$x_1,x_2\in B^2(z_m,a_m)$, then the inequality
$|x_1-z_m|\geq|x_2-z_m|$ implies
$$
m_2(B^2(z_m,a_m)\cap B^2(x_1,r))\leq m_2(B^2(z_m,a_m)\cap
B^2(x_2,r)).
$$
Thus we have
$$
C=\inf_{r\leq a_m}\frac{m_2(B^2(z_m,a_m)\cap
B^2(z_m+a_m,r))}{m_2(B^2(z_m+a_m,r))}.
$$
The right-hand side of the last formula is invariant under the
similarities. Consequently using the similarity
$$
\mathbb C\ni z\longmapsto\frac1r(z-z_m)\in\mathbb C
$$
we see that
\begin{multline}\label{20*}
C=\inf_{r\leq a_m}\frac{m_2(B^2(0,\frac{a_m}r)\cap
B^2(\frac{a_m}r,1))}{m_2(B^2(\frac{a_m}r,1))}=\inf_{r\geq1}
\frac{m_2(B^2(0,r)\cap
B^2(r,1))}{m_2(B^2(r,1))}\\=\frac1\pi\inf_{r\geq1}m_2(B^2(0,r)\cap
B^2(r,1)) =\frac 1\pi m_2(B^2(-1,1)\cap
B^2(0,1))=\frac23-\frac{\sqrt{3}}{2\pi}.
\end{multline}
\begin{figure}[tbh]
\includegraphics{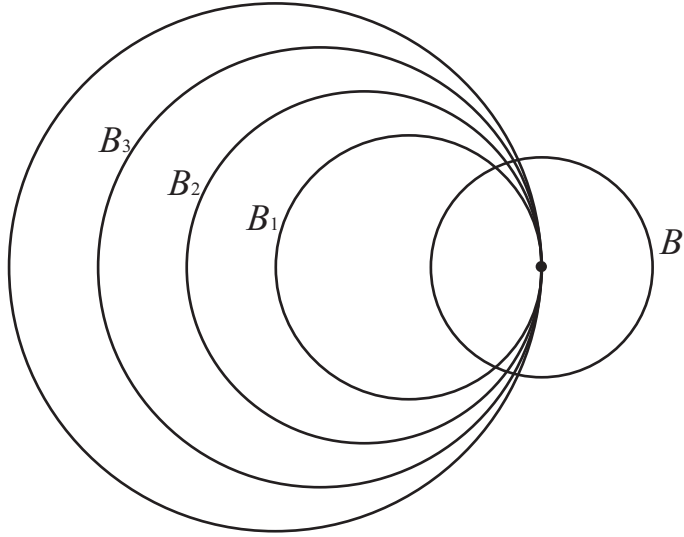}
\caption{The center of the fixed small ball $B$ lies on the
boundary spheres of the large balls $B_n$. The volume of the
intersection $B\cap B_n$  grows together with the radius of the
ball $B_n$.}
\end{figure}
The last equality, \eqref{e18} and \eqref{e19} imply that
$$
\frac1{(\frac23-\frac{\sqrt{3}}{2\pi})m_2(B^2(x,r))}\int\limits_{B^2(x,r)}u_3(y)\,
dm_2(y)\geq u_3(x)
$$
whenever $r\in A$ and $\overline{B^2(x,r)}\subseteq\Omega$.

Thus the theorem is proved in the case where limit relations
\eqref{eq12} and \eqref{eq13} hold. Similar constructions can be
realized if \eqref{eq12} and \eqref{eq14} hold and we omit them
here.
\end{proof}

Statement (s) of Theorem \ref{t3.6} has a useful reformulation.
For $A\subseteq(0,\infty)$ define
$$
\ln(A):=\{\ln x:x\in A\}
$$
with $\ln(\emptyset):=\emptyset$. Then $\ln(A)$ is a subset of
$\mathbb R$. Recall that a set $X\subseteq\mathbb R$ is an
$\varepsilon$-net in $\mathbb R,\ \varepsilon>0$, if
$$
\mathbb R=\bigcup_{x\in X}B^1(x,\varepsilon).
$$

\subsection{Proposition}\label{p2.9}{\it Let $A$ be a subset of $(0,\infty)$.
Then statement (s) in Theorem \ref{t3.6} is valid with this $A$ if
and only if there is $\varepsilon>0$ such that $\ln(A)$ is an
$\varepsilon$-net in $\mathbb R$.}
\begin{proof}
If (s) holds, then $\ln(A)$ is an $\varepsilon$-net with
$\varepsilon=\ln C$ where $C$ is the constant in \eqref{e9}. If
(s) does not hold, then Lemma \ref{l3.8} implies that $\ln A$ is
not an $\varepsilon$-net for any $\varepsilon>0$.
\end{proof}
Using this proposition and analysing the first part of the proof
of Theorem \ref{t3.6} we obtain the following
\subsection{Proposition}\label{pppp}
{\it Let $A$ be a subset of $(0, \infty)$. The following three
statements are equivalent.

\begin{itemize}
\item[(i)] $A$ is favorable for all domains of ${\mathbb R}^n$.

\item[(ii)] $A$ is favorable for all open sets of  ${\mathbb
R}^n$.

\item[(iii)] There is $\varepsilon>0$ such that $\ln(A)$ is an
$\varepsilon$-net in ${\mathbb R}$.
\end{itemize}}

The condition for the set $A\subseteq(0,\infty)$, to be favorable
for all {\it bounded} domains $\Omega$ can be presented in terms
of porosity of $A$, so recall a definition.
\subsection{Definition}\label{d2.10}{\it Let $A\subseteq(0,\infty)$.
The right hand porosity of $A$ at zero is the quantity $$
p_0(A):=\limsup_{h\to0+}\frac{l(h,A)}h
$$
where $l(h,A)$ is the length of the longest interval in
$[0,h]\setminus A,\ h>0$.}
\subsubsection{Remark}\label{r2.10.1}{\rm
It is easy to see that $0\leq p_0(A)\leq1$ for each
$A\subseteq\mathbb R$. A variety of computations directly related
to the notion of porosity can be found in \cite[p.~183--212]{Th}.}

\subsection{Definition}\label{d2.10*}{\it
Let $A\subseteq(0,\infty)$. The right porosity index of $A$ at
$0$, $i_0(A)$, is defined to be the supremum of all real numbers
$r$ for which there is a sequence of open intervals
$\{(a_n,b_n)\}_{n\in\mathbb N},\ a_n<b_n$, such that
$\lim_{n\to\infty}a_n=\lim_{n\to\infty}b_n=0$ and
$(a_n,b_n)\subset(0,\infty)\setminus A$ and
$$
r<\frac{b_n-a_n}{a_n}
$$
for each $n\in\mathbb N$. }

 If no such numbers $r$ exist, then following the usual
 conversion we define $i_0(A):=0$.

 The following lemma is a particular case of Lemma $A_{2.13}$ from
 \cite[p.~185]{Th}.
 \subsection{Lemma}\label{l2.11*}{\it
The equality
$$
i_0(A)=\frac{p_0(A)}{1-p_0(A)}
$$
holds for each $A\subseteq(0,\infty)$. }

\subsection{Theorem}\label{t2.11}{\it Let $A$ be a subset of $(0,\infty)$.
Then $A$ is favorable for all bounded domains
$\Omega\subseteq\mathbb R^n$ if and only if $ p_0(A)<1. $}

\begin{proof}
It follows from Lemma \ref{l2.11*} that $p_0(A)=1$ if and only if
$i_0(A)=1$. Using the definition of porosity index $i_0(A)$ we can
prove that the equality $i_0(A)=\infty$ implies the existence of
disjoint intervals $(a_n,b_n)\subset(0,\infty)\setminus A,\
n=1,2,\dots,$ such that equations \eqref{eq12} and \eqref{eq13}
hold. It was shown in the proof of Theorem~\ref{t3.6} that if
\eqref{eq12} and \eqref{eq13} hold then
 there are a domain
$\Omega\subseteq\mathbb R^n$ and a nonnegative $u\in\mathcal
L_{\rm loc}^1(\Omega)\setminus QNS(\Omega)$ such that \eqref{e8}
holds whenever $r\in A$ and $\overline{B^n(x,r)}\subseteq\Omega$.
It remains to observe that inequality \eqref{e13} implies
$diam(\Omega)<\infty$. Thus if $A$ is favorable for all bounded
domains $\Omega\subseteq\mathbb R^n$, then $p_0(A)<1$.

Now note that if $p_0(A)<1$, then the set $(-\infty,R)\cap\ln(A)$
is an $\varepsilon$-net, $\varepsilon=\varepsilon(R)$, in
$(-\infty,R)$ for each $R\in\mathbb R$. Hence, reasoning as in the
first part of the proof of Theorem \ref{t3.6} we can prove the
implication $$(p_0(A)<1)\Rightarrow(A \text{ is favorable for
every bounded domain }\Omega\subseteq\mathbb R^n).
$$
\end{proof}

\subsubsection{Remark}\label{r2.11.1} {\rm As in Proposition \ref{pppp} it is easy
to prove that $A$ is favorable for all bounded domains of $\mathbb
R^n$ if and only if $A$ is favorable for all bounded open subsets
of $\mathbb R^n$. Theorem \ref{t2.11} remains valid even for
unbounded domains and open sets $\Omega\subseteq\mathbb R^n$ if
$$
\sup_{x\in\Omega}(\delta_\Omega(x))<\infty.
$$}

\subsubsection{Remark}\label{r2.11.2} {\rm In complete analogy with
Definition \ref{d2.10*} we may define the quantity $i_\infty(A)$,
the left porosity index of $A$ at $+\infty$, after which
Proposition~\ref{p2.9} can be reformulated as:}

{\it Let $A\subset(0,\infty)$. The following statements are
equivalent:
\begin{itemize}
\item[(i)] $A$ is favorable for all domains
$\Omega\subseteq\mathbb R^n$;

\item[(ii)] The indexes $i_0(A)$ and $i_\infty(A)$ are less than
infinity,
$$
i_0(A)\vee i_\infty(A)<\infty;
$$

\item[(iii)] There is $\varepsilon>0$ such that $\ln(A)$ is an
$\varepsilon$-net in $\mathbb R^1$.
\end{itemize} }
\medskip

Theorem \ref{t3.6}, Proposition \ref{p2.9} and Theorem \ref{t2.11}
imply the following.

\subsection{Corollary}\label{c2.12}{\it Let $A$ be a subset of $(0,\infty)$
and let $\alpha,\beta$ be positive constants. Then the set $A$ is
favorable for all domains $\Omega\subseteq\mathbb R^n$ (for all
bounded domains $\Omega\subseteq\mathbb R^n$) if and only if the
set $$ \alpha A^{\beta}:=\{\alpha x^\beta:x\in A\}
$$ has the same property.}
\begin{proof}
 One just directly observe that if   condition (s) holds for the
 set $A$ with a constant $C$, then   condition (s) holds for the
 set $\alpha A^\beta$ with the constant $C'\geq C^\beta$.
\end{proof}
{\rm Now we are ready to characterize the function
$f:(0,\infty)\to(0,\infty)$ for which the equality
\begin{equation}\label{e21}
Q(f,D,\Omega)=QNS(\Omega)
\end{equation}
is fulfilled for all open sets $\Omega\subseteq\mathbb R^n$.}

\subsection{Theorem}\label{p2.13}{\it Let $D$ be a bounded, Lebesgue measurable
subset of $\mathbb R^n$ with a marked point $p_D\in Int\,D$ and
let $\varphi:SM(\mathbb R^n) \to (0,\infty)$,
$f:(0,\infty)\to(0,\infty)$ be functions such that \eqref{e5}
holds for all $h\in SM(\mathbb R^n)$. Then  equality \eqref{e21}
holds for all open sets $\Omega\subseteq\mathbb R^n$ if and only
if there are $A\subseteq(0,\infty)$ and $c>1$ such that:
\begin{itemize}

\item[(i)] the inequality $f(k)\leq ck$  holds for all
$k\in(0,\infty)$,

\item[(ii)] $\ln(A)$ is an $\varepsilon$-net in $\mathbb R$ for
some $\varepsilon>0$,

\item[(iii)] the inequality
\begin{equation}\label{e22}
\frac1ck\leq f(k)
\end{equation}
holds for all $k\in A$.
\end{itemize}}
\begin{proof}
Let $\Omega$ be an open set in $\mathbb R^n$,
$A\subseteq(0,\infty)$ and $c>1$. Assume that $\ln(A)$ is an
$\varepsilon$-net in $\mathbb R$ for some $\varepsilon>0$ and that
\eqref{e22} holds with this $c$ for all $k\in A$. Then using
\eqref{e6} and \eqref{e22} we obtain
\begin{equation}\label{e23}
u(x_\Omega)\leq\frac{c^n K(u)}{(k(h))^n}\int\limits_{h(D)}
u(y)\,dm_n(y)
\end{equation}
for every $u\in Q(f,D,\Omega)$ and every $x_\Omega$ whenever $h\in
Sim(p_D,x_\Omega)$ and $k(h)\in A$. As in the proof of Theorem~2.2
write
$$
R_D=\sup_{y\in D}|p_D-y|.
$$
Let $B^n(x_\Omega,r_0)$ be a ball such that $r_0=R_Dk_0$, with
$k_0\in A$ and $\overline{B^n(x_\Omega,r_0)}\subset\Omega$. Then
each similarity $h$ such that $k(h)=k_0$ and  $h(p_D)=x_\Omega$
belongs to $Sim(p_D,x_\Omega)$ and satisfies $h(D)\subseteq
B^n(x_\Omega,r_0)$. Consequently \eqref{e23} implies
$$
u(x_\Omega)\leq\frac{c^n
K(u)(R_D)^n\nu_n}{(k_0R_D)^n\nu_n}\int\limits_{B^n(x_\Omega,
r_0)}u(y)\,dm_n(y)= \frac{c^n K(u)\nu_n(R_D)^n}{
m_n(B^n(x_\Omega,r_0))}
\int\limits_{B^n(x_\Omega,r_0)}u(y)\,dm_n(y)
$$
for every $u\in Q(f,D,\Omega)$ and every
$B^n(x_\Omega,r_0)\subseteq\Omega$ whenever $r_0\in R_DA$.
Corollary \ref{c2.12} implies that the set $R_DA$ is favorable for
$\Omega$. Hence $Q(f,D,\Omega)\subseteq QNS(\Omega)$. Taking into
account Proposition 2.4 we see that  conditions (i)--(iii) of the
present theorem imply equality \eqref{e21} for all open sets
$\Omega\subseteq\mathbb R^n$.

Conversely, suppose that \eqref{e21} holds for all open sets
$\Omega\subseteq\mathbb R^n$ but for every $t>0$ the set $\ln
A_t$, where
\begin{equation}\label{e24}
A_t:=\{k\in(0,\infty):f(k)\geq tk\},
\end{equation}
is not an $\varepsilon$-net for any $\varepsilon>0$. It is clear
that $A_{t_1}\subseteq A_{t_2}$ if $t_1\geq t_2$. Applying
Proposition \ref{p2.9} and Lemma \ref{l3.8} to the sets
$A_2,A_3,\dots$ we  obtain a sequence $\{a_m\}_{m=2}^\infty$ of
positive numbers $a_m$ such that $A_m\cap(a_m,ma_m)=\emptyset$ for
each $m\geq2$, i.e.,
\begin{equation}\label{e25}
f(k)<\frac1mk
\end{equation}
if $a_m<k<ma_m$, and that
\begin{equation}\label{e26}
(a_{m_1},m_1a_{m_1})\cap(a_{m_2},m_2a_{m_2})=\emptyset
\end{equation}
whenever $m_1\ne m_2$. Passing, if necessary, to a subsequence we
may assume that $\{a_m\}_{m=2}^\infty$ and $\{ma_m\}_{m=2}^\infty$
are monotone and convergent in $[0,\infty]$ sequences. This
assumption and \eqref{e26} imply either the equalities
\begin{equation}\label{e27}
\lim_{m\to\infty}a_m=\lim_{m\to\infty}ma_m=0
\end{equation}
or the equalities
$$
\lim_{m\to\infty}a_m=\lim_{m\to\infty}ma_m=\infty.
$$
As in the proof of Theorem \ref{t3.6} we consider only the case
when \eqref{e27} holds and  the dimension $n=2$. We shall
construct a domain $\Omega\subseteq\mathbb R^2$ and a nonnegative
$u\in\mathcal L_{\rm loc}^1(\Omega)$ such that
$$
u\in Q(f,D,\Omega)\setminus QNS(\Omega).
$$
To this end, note that \eqref{e21} implies \eqref{e7}, so using
Proposition 2.4 we can find $c\geq1$ such that
\begin{equation}\label{e28}
f(k)\leq ck
\end{equation}
for every $k\in(0,\infty)$. Let us define a function
$f_1:(0,\infty)\to(0,\infty)$ by the rule
\begin{equation}\label{e29}
f_1(k):=\begin{cases} \frac km&\text{if }a_m<k<ma_m,\ m=2,3,\dots\\
ck&\text{if }k\in(0,\infty)\setminus\bigcup_{m=2}^\infty(a_m,ma_m)
\end{cases}
\end{equation}
where $c\geq1$ is the constant from inequality \eqref{e28}.
Inequalities \eqref{e25} and \eqref{e28} imply $f(k)\leq f_1(k)$
for all $k\in(0,\infty)$. Hence from the definition of the set
$Q(f,D,\Omega)$ follows the inclusion
$$
Q(f,D,\Omega)\supseteq Q(f_1,D,\Omega).
$$
Thus it is sufficient to find a domain $\Omega\subseteq\mathbb
R^2$ and a nonnegative $u\in\mathcal L_{\rm loc}^1(\Omega)$ such
that
$$
u\in Q(f_1,D,\Omega)\setminus QNS(\Omega).
$$
Let us define
$$
\Omega :=\bigcup_{m=N_0+1}^\infty \left(B^2\left(z_m,
\frac{ma_m}{N_0}\right)\cup R_m\right), \ \ \ u(x):=\left\{
\begin{aligned}
&1\ \ \mbox{if}\ \ \ x\in X\\
&0\ \ \mbox{if}\ \ \ x\in \Omega\setminus X
\end{aligned}
\right.
$$
where
$$
X:=\bigcup_{m=N_0+1}^\infty \overline{B^2(z_m, a_m)},\ \ \
z_m:=2\sum_{i=1}^{m-1} i a_i
$$
$$
R_m:=\{z\in {\mathbb C} : 2\sum_{n=1}^{m-1} n a_n < Re\,
(z)<2\sum_{n=1}^{m} n a_n, \ \ \ |Im(z)|< a_{m+1}\}.
$$
The parameter $N_0$ is free here and we will specify this
parameter later. It is relevant to remark that the domain $\Omega$
is obtained from the domain depicted on Fig.~1, by deleting of the
balls $B^2(z_1,\frac{b_1}{N_0}),
B^2(z_2,\frac{b_2}{N_0}),\dots,\break
B^2(z_{N_0},\frac{b_{N_0}}{N_0})$ and the rectangels $R_1,\dots,
R_{N_0}$ and putting $b_m:=ma_m$ in the rest of balls and
rectangels. As in the proof of Theorem~\ref{t3.6} we have
$u\not\in QNS(\Omega)$. It still remains to prove that $u\in
Q(f_1,D,\Omega)$. The last relation holds if and only if there
exists $K(u)\geq1$ such that
\begin{equation}\label{e30}
(f_1(k(h)))^2\leq K(u)\int\limits_{h(D)}u(y)\,dm_2(y)
\end{equation}
for   all $h\in Sim(p_D,x_\Omega)$ with $x_\Omega\in X$.

Let $x_\Omega\in X$. It follows from the definitions of $\Omega$
and $X$ that there is $m\geq N_0+1$ for which
\begin{equation*}\label{e31}
x_\Omega\in\overline{B^2(z_m,a_m)}.
\end{equation*}

We claim that the inequality
\begin{equation}\label{e32}
k(h)r_D\leq\frac{2ma_m}{N_0}.
\end{equation}
holds for every $h\in Sim(p_D,x_\Omega)$ with
$r_D=\delta_{Int(D)}(p_D)$.

Let us prove it. Since $h\in Sim(p_D,x_\Omega)$ we have
$$
h(B^2(p_D,r_D))\subseteq \Omega.
$$
The last inclusion implies
$$
\partial\Omega\cap h(B^2(p_D,r_D))=\emptyset
$$
because $\Omega\cap\partial \Omega=\emptyset$ for the open sets.
The intersection
$$
\partial B^2\Big(z_m,\frac{ma_m}{N_0}\Big)
\cap\partial\Omega=\Big\{z\in\mathbb C:|z-z_m|=\frac{ma_m}{N_0}
\Big\}\cap\partial\Omega
$$
is not empty, see Fig.~1. Consequently there is $\xi\in \partial
B^2(z_m,\frac{ma_m}{N_0} )\setminus h(B^2(p_D,r_D))$. Hence
$$
|x_\Omega-\xi|\geq k(h)r_D.
$$
Using the triangle inequality we obtain
$$
|x_\Omega-\xi|\leq
|x_\Omega-z_m|+|z_m-\xi|=|x_\Omega-z_m|+\frac{ma_m}{N_0}.
$$
Consequently
$$
k(h)r_D\leq|x_\Omega-z_m|+\frac{ma_m}{N_0}.
$$
Since $x_\Omega\in\overline{B^2(z_m,a_m)}$ we have
$|x_\Omega-z_m|\leq a_m$. It follows directly from the definition
of $\Omega$ that $m\geq N_0$. Hence
$$
k(h)r_D\leq a_m+\frac{ma_m}{N_0}\leq \frac{2ma_m}{N_0}.
$$
Inequality \eqref{e32} follows.

Since $h(D)\supseteq h(B^2(p_D,r_D))$, the inequality
\begin{equation}\label{e33}
(f_1(k(h)))^2\leq
K(u)\int\limits_{B^2(x_\Omega,k(h)r_D)}u(y)\,dm_2(y)
\end{equation}
implies \eqref{e30}, so it is sufficient to prove \eqref{e33}. The
following two cases are possible: $k(h)\in(0,a_m]$ and
$k(h)\in(a_m,\infty)$. Before analyzing these cases note that $
f_1(k)\leq ck$ for every $k\in (0, \infty)$ because $\frac 1m\leq
\frac 12$ and $c\geq1$ in definition \eqref{e29}. Hence in
 the first
case  we can replace \eqref{e33} by
\begin{equation}\label{e34}
c^2\leq\frac{K(u)}{(k(h))^2}\int\limits_{B^2(x_\Omega,k(h)r_D)}u(y)\,dm_2(y).
\end{equation}
It is clear that
\begin{equation}\label{e35}
\frac{K(u)}{(k(h))^2}\int\limits_{B^2(x_\Omega,k(h)r_D)}u(y)\,dm_2(y)\geq
\frac{K(u)\pi(r_D\wedge1)^2}{\pi(k(h)(r_D\wedge1))^2}\int\limits_
{B^2(x_\Omega,k(h)(r_D\wedge1))}u(y)\,dm_2(y).
\end{equation}
Since $k(h)\in(0,a_m]$, we see that
$$
k(h)(r_D\wedge1)\leq a_m.
$$
Hence, as it was shown in the proof of Theorem \ref{t3.6}, in the
case under consideration we have
\begin{equation*}
\frac1{\pi(k(h)(r_D\wedge1))^2}\int\limits_
{B^2(x_\Omega,k(h)(r_D)\wedge1)}u(y)\,dm_2(y)\geq\frac23-\frac{\sqrt{3}}{2\pi}.
\end{equation*}
The last estimation and \eqref{e35} show that \eqref{e34} holds if
$$
c^2=K(u)\pi(r_D\wedge1)^2\Big(\frac23-\frac{\sqrt{3}}{2\pi}\Big).
$$
Consider now the case $k(h)\in(a_m,\infty)$. Inequality
\eqref{e32} shows that
$$
k(h)\leq\frac{2ma_m}{N_0r_D}.
$$
Let us specify $N_0$ as the smallest positive integer $N$
satisfying the inequality $\frac2{Nr_D}>1$. Then we obtain the
double inequality
$$
a_m<k(h)<ma_m.
$$
This inequality and \eqref{e29} show that
$$
f_1(k(h))=\frac km\leq a_m.
$$
Consequently we can prove \eqref{e33} as in the case
$k(h)\in(0,a_m]$.
\end{proof}

Let us consider now some examples of functions $\varphi $ and $f$
for which equality \eqref{e5} holds.

\subsection{Example}\label{e2.15} Let $\psi$ be a positive bounded
periodic function on $\mathbb R$. Write $$
\mu(x):=\frac12\Big(x+\frac1x\Big)
$$ for $x>0$ and  define
\begin{equation}\label{e40}
f(k):=k\psi(\mu(k)).
\end{equation}
Using some routine estimations we see that conditions (i)--(iii)
from Theorem \ref{p2.13} are satisfied by the function $f$ if we
take
$$
A=\mu^{-1}\{x\in(0,\infty):\psi(x)\geq\frac12M\}, \quad
c=M\vee\frac2M
$$
where
$$
M=\sup_{y\in\mathbb R}\psi(y).
$$

An important special case of the preceding example is the constant
function $\psi$. Then $f$ is   linear  on $(0,\infty)$ and
conditions (i)--(iii) from Theorem \ref{p2.13} are evidently hold.
In this simplest case the  function $\varphi:SM(\mathbb
R^n)\to(0,\infty)$ connected with $f$ can be obtained in distinct
ways depending on the geometrical properties of the set $D$.

In all following examples $D$ is a bounded Lebesgue measurable
subset of $\mathbb R^n$ with $Int\,D\ne\emptyset$ and $h\in
SM(\mathbb R^n)$.

\subsection{Example}\label{e2.17.2} Let $d$-dimensional Hausdorff
measure $\mathcal H^d$, $n-1\leq d\leq n$, of the boundary
$\partial D$ be finite and nonzero, $0<\mathcal H^d(\partial
D)<\infty$. Write
$$
\varphi(h)=(\mathcal H^d(\partial(h(D))))^{\frac1d}.
$$

\subsection{Example}\label{e2.17.3} Let $D$ be a set with the
finite Caccoppoli--de Gorgi perimeter $P$, see, for instance,
\cite[Chapter~3, \S3]{BZ}. Write
$$
\varphi(h)=(P(h(D)))^{\frac1{n-1}}.
$$

\subsection{Example}\label{e2.17.4} Let $D\subseteq\mathbb R^2$ be a simply connected
domain with  the rectifiable boundary $\partial D,\ 0<\mathcal
H^1(\partial D)<\infty$. Suppose that the domain $D$ is not a
disk. Write
$$
\varphi(h)=((\mathcal H^1(h(\partial D)))^2-4\pi
m_2(h(D)))^{\frac12}.
$$
In this case the inequality $\varphi(h)>0$ follows from the
Classical Isoperemetric Inequality, see, for instance,
\cite[Chapter~1,\S1]{BZ}.

\medskip

This list of examples can be simply extended by involving the
analytic capacity, the transfinite diameter, the Menger curvature
etc. for the definition of the function $\varphi$. The homogenity
under dilatations $x\longmapsto\alpha x,\ x\in\mathbb R^n,\
\alpha>0$, the invariance under isometries, finiteness and
positiveness are sufficient for this purpose.

\medskip

{\bf Acknowledgment.} The first author was partially supported by
the Academy of Finland.

\end{document}